\input phyzzx
\catcode`@=11 
\def\space@ver#1{\let\@sf=\empty \ifmmode #1\else \ifhmode
   \edef\@sf{\spacefactor=\the\spacefactor}\unskip${}#1$\relax\fi\fi}
\def\attach#1{\space@ver{\strut^{\mkern 2mu #1} }\@sf\ }
\newtoks\foottokens
\newbox\leftpage \newdimen\fullhsize \newdimen\hstitle \newdimen\hsbody
\newif\ifreduce  \reducefalse
\def\almostshipout#1{\if L\lr \count2=1
      \global\setbox\leftpage=#1 \global\let\lr=R
  \else \count2=2
    \shipout\vbox{\special{dvitops: landscape}
      \hbox to\fullhsize{\box\leftpage\hfil#1}} \global\let\lr=L\fi}
\def\smallsize{\relax
\font\eightrm=cmr8 \font\eightbf=cmbx8 \font\eighti=cmmi8
\font\eightsy=cmsy8 \font\eightsl=cmsl8 \font\eightit=cmti8
\font\eightt=cmtt8
\def\eightpoint{\relax
\textfont0=\eightrm  \scriptfont0=\sixrm
\scriptscriptfont0=\sixrm
\def\rm{\fam0 \eightrm \f@ntkey=0}\relax
\textfont1=\eighti  \scriptfont1=\sixi
\scriptscriptfont1=\sixi
\def\oldstyle{\fam1 \eighti \f@ntkey=1}\relax
\textfont2=\eightsy  \scriptfont2=\sixsy
\scriptscriptfont2=\sixsy
\textfont3=\tenex  \scriptfont3=\tenex
\scriptscriptfont3=\tenex
\def\it{\fam\itfam \eightit \f@ntkey=4 }\textfont\itfam=\eightit
\def\sl{\fam\slfam \eightsl \f@ntkey=5 }\textfont\slfam=\eightsl
\def\bf{\fam\bffam \eightbf \f@ntkey=6 }\textfont\bffam=\eightbf
\scriptfont\bffam=\sixbf   \scriptscriptfont\bffam=\sixbf
\def\tt{\fam\ttfam \eightt \f@ntkey=7 }
\def\caps{\fam\cpfam \tencp \f@ntkey=8 }\textfont\cpfam=\tencp
\setbox\strutbox=\hbox{\vrule height 7.35pt depth 3.02pt width\z@}
\samef@nt}
\def\Eightpoint{\eightpoint \relax
  \ifsingl@\subspaces@t2:5;\else\subspaces@t3:5;\fi
  \ifdoubl@ \multiply\baselineskip by 5
            \divide\baselineskip by 4\fi }
\parindent=16.67pt
\itemsize=25pt
\thinmuskip=2.5mu
\medmuskip=3.33mu plus 1.67mu minus 3.33mu
\thickmuskip=4.17mu plus 4.17mu
\def\thinspace{\kern .13889em }
\def\negthinspace{\kern-.13889em }
\def\enspace{\kern.416667em }
\def\enskip{\hskip.416667em\relax}
\def\quad{\hskip.83333em\relax}
\def\qquad{\hskip1.66667em\relax}
\def\crr{\cropen{8.3333pt}}
\foottokens={\Eightpoint\singlespace}
\def\papersize{\SIZE\OFFSET\skip\footins=\bigskipamount}
\def\SIZE{\hsize=11.8truecm\vsize=17.5truecm}
\def\OFFSET{\voffset=-1.3truecm\hoffset=  .14truecm}
\message{STANDARD CERN-PREPRINT FORMAT}
\def\attach##1{\space@ver{\strut^{\mkern 1.6667mu ##1} }\@sf\ }
\def\PH@SR@V{\doubl@true\baselineskip=20.08pt plus .1667pt minus .0833pt
             \parskip = 2.5pt plus 1.6667pt minus .8333pt }
\def\author##1{\vskip\frontpageskip\titlestyle{\tencp ##1}\nobreak}
\def\address##1{\par\kern 4.16667pt\titlestyle{\tenpoint\it ##1}}
\def\andaddress{\par\kern 4.16667pt \centerline{\sl and} \address}
\def\abstract{\vskip\frontpageskip\centerline{\twelverm ABSTRACT}
              \vskip\headskip }
\def\cases##1{\left\{\,\vcenter{\Tenpoint\m@th
    \ialign{$####\hfil$&\quad####\hfil\crcr##1\crcr}}\right.}
\def\matrix##1{\,\vcenter{\Tenpoint\m@th
    \ialign{\hfil$####$\hfil&&\quad\hfil$####$\hfil\crcr
      \mathstrut\crcr\noalign{\kern-\baselineskip}
     ##1\crcr\mathstrut\crcr\noalign{\kern-\baselineskip}}}\,}
\Tenpoint
}
\def\Smallsize{\smallsize\reducetrue
\let\lr=L
\hstitle=8truein\hsbody=4.75truein\fullhsize=24.6truecm\hsize=\hsbody
\output={
  \almostshipout{\leftline{\vbox{\makeheadline
  \pagebody\makefootline}}}\advancepageno
     }
\special{dvitops: landscape}
\def\makeheadline{
\iffrontpage\line{\the\headline}
             \else\vskip .0truecm\line{\the\headline}\vskip .5truecm \fi}
\def\makefootline{\iffrontpage\vskip  0.truecm\line{\the\footline}
               \vskip -.15truecm\line{\the\date\hfil}
              \else\line{\the\footline}\fi}
\paperheadline={
\iffrontpage\hfil
               \else
               \tenrm\hss $-$\ \folio\ $-$\hss\fi    }
\paperstyle}
\newcount\referencecount     \referencecount=0
\newif\ifreferenceopen       \newwrite\referencewrite
\newtoks\rw@toks
\def\NPrefmark#1{\attach{\scriptscriptstyle [ #1 ] }}
\let\PRrefmark=\attach
\def\refmark#1{\relax\ifPhysRev\PRrefmark{#1}\else\NPrefmark{#1}\fi}
\def\refend{\refmark{\number\referencecount}}
\newcount\lastrefsbegincount \lastrefsbegincount=0
\def\refsend{\refmark{\count255=\referencecount
   \advance\count255 by-\lastrefsbegincount
   \ifcase\count255 \number\referencecount
   \or \number\lastrefsbegincount,\number\referencecount
   \else \number\lastrefsbegincount-\number\referencecount \fi}}
\def\refch@ck{\chardef\rw@write=\referencewrite
   \ifreferenceopen \else \referenceopentrue
   \immediate\openout\referencewrite=referenc.texauxil \fi}
%
{\catcode`\^^M=\active 
  \gdef\obeyendofline{\catcode`\^^M\active \let^^M\ }}%
%
{\catcode`\^^M=\active 
  \gdef\ignoreendofline{\catcode`\^^M=5}}
{\obeyendofline\gdef\rw@start#1{\def\t@st{#1} \ifx\t@st\blankend%
\endgroup \@sf \relax \else \ifx\t@st\bl@nkend \endgroup \@sf \relax%
\else \rw@begin#1
\backtotext
\fi \fi } }
{\obeyendofline\gdef\rw@begin#1
{\def\n@xt{#1}\rw@toks={#1}\relax%
\rw@next}}
\def\blankend{}
{\obeylines\gdef\bl@nkend{
}}
\newif\iffirstrefline  \firstreflinetrue
\def\rwr@teswitch{\ifx\n@xt\blankend \let\n@xt=\rw@begin %
 \else\iffirstrefline \global\firstreflinefalse%
\immediate\write\rw@write{\noexpand\obeyendofline \the\rw@toks}%
\let\n@xt=\rw@begin%
      \else\ifx\n@xt\rw@@d \def\n@xt{\immediate\write\rw@write{%
        \noexpand\ignoreendofline}\endgroup \@sf}%
             \else \immediate\write\rw@write{\the\rw@toks}%
             \let\n@xt=\rw@begin\fi\fi \fi}
\def\rw@next{\rwr@teswitch\n@xt}
\def\rw@@d{\backtotext} \let\rw@end=\relax
\let\backtotext=\relax

\newdimen\refindent     \refindent=30pt
\def\refitem#1{\par \hangafter=0 \hangindent=\refindent \Textindent{#1}}
\def\REFNUM#1{\space@ver{}\refch@ck \firstreflinetrue%
 \global\advance\referencecount by 1 \xdef#1{\the\referencecount}}
\def\refnum#1{\space@ver{}\refch@ck \firstreflinetrue%
 \global\advance\referencecount by 1 \xdef#1{\the\referencecount}\refend}

\def\REF#1{\REFNUM#1%
 \immediate\write\referencewrite{%
 \noexpand\refitem{#1.}}%
\begingroup\obeyendofline\rw@start}
\def\ref{\refnum\?%
 \immediate\write\referencewrite{\noexpand\refitem{\?.}}%
\begingroup\obeyendofline\rw@start}
\def\Ref#1{\refnum#1%
 \immediate\write\referencewrite{\noexpand\refitem{#1.}}%
\begingroup\obeyendofline\rw@start}
\def\REFS#1{\REFNUM#1\global\lastrefsbegincount=\referencecount
\immediate\write\referencewrite{\noexpand\refitem{#1.}}%
\begingroup\obeyendofline\rw@start}
\newif\ifmref  
\newif\iffref  
\def\xrefsend{\xrefmark{\count255=\referencecount
\advance\count255 by-\lastrefsbegincount
\ifcase\count255 \number\referencecount
\or \number\lastrefsbegincount,\number\referencecount
\else \number\lastrefsbegincount-\number\referencecount \fi}}
\def\xrefsdub{\xrefmark{\count255=\referencecount
\advance\count255 by-\lastrefsbegincount
\ifcase\count255 \number\referencecount
\or \number\lastrefsbegincount,\number\referencecount
\else \number\lastrefsbegincount,\number\referencecount \fi}}
\def\xREFNUM#1{\space@ver{}\refch@ck\firstreflinetrue%
\global\advance\referencecount by 1
\xdef#1{\xrefend}}
\def\xrefend{\xrefmark{\number\referencecount}}
\def\xrefmark#1{[{#1}]}
\def\xRef#1{\xREFNUM#1\immediate\write\referencewrite%
{\noexpand\refitem{#1 }}\begingroup\obeyendofline\rw@start}%
\def\xREFS#1{\xREFNUM#1\global\lastrefsbegincount=\referencecount%
\immediate\write\referencewrite{\noexpand\refitem{#1 }}%
\begingroup\obeyendofline\rw@start}
\def\rrr#1#2{\relax\ifmref{\iffref\xREFS#1{#2}%
\else\xRef#1{#2}\fi}\else\xRef#1{#2}\xrefend\fi}
\def\multref#1#2{\mreftrue\freftrue{#1}%
\freffalse{#2}\mreffalse\xrefsend}
\referencecount=0
\def\refbreak{\hfil\penalty200\hfilneg}
\def\paperstyle{\papers}
\paperstyle   
%
%
%
\def\slacpub{\afterassignment\slacp@b\toks@}
\def\slacp@b{\edef\n@xt{\Pubnum={NIKHEF/\the\toks@}}\n@xt}
\let\pubnum=\slacpub
\expandafter\ifx\csname eightrm\endcsname\relax
    \let\eightrm=\ninerm \let\eightbf=\ninebf \fi

\font\seventeencp=cmcsc10 scaled\magstep3

\newif\ifCONF \CONFfalse
\newif\ifBREAK \BREAKfalse
\newif\ifsectionskip \sectionskiptrue

%
%
%
%
\def\NuclPhysProc{
\let\lr=L
\hstitle=8truein\hsbody=4.75truein\fullhsize=21.5truecm\hsize=\hsbody
\hstitle=8truein\hsbody=4.75truein\fullhsize=20.7truecm\hsize=\hsbody
\output={
  \almostshipout{\leftline{\vbox{\makeheadline
  \pagebody\makefootline}}}\advancepageno
     }
\def\papersize{\SIZE\OFFSET\skip\footins=\bigskipamount}
\def\SIZE{\hsize=10.0truecm\vsize=27.0truecm}
\def\OFFSET{\voffset=-1.4truecm\hoffset=-2.40truecm}
\message{NUCLEAR PHYSICS PROCEEDINGS FORMAT}
\def\makeheadline{
\iffrontpage\line{\the\headline}
             \else\vskip .0truecm\line{\the\headline}\vskip .5truecm \fi}
\def\makefootline{\iffrontpage\vskip  0.truecm\line{\the\footline}
               \vskip -.15truecm\line{\the\date\hfil}
              \else\line{\the\footline}\fi}
\paperheadline={\hfil}
\paperstyle}
%
%
%
%

%
%
%
%

%
%
%
%
\def\ReprintVolume{\smallsize
\def\papersize{\hsize=18.0truecm\vsize=25.1truecm\voffset -.73truecm
    \hoffset -.65truecm\skip\footins=\bigskipamount
    \normaldisplayskip= 20pt plus 5pt minus 10pt}
\message{REPRINT VOLUME FORMAT}
\paperstyle\baselineskip=.425truecm\parskip=0truecm
\def\makeheadline{
\iffrontpage\line{\the\headline}
             \else\vskip .0truecm\line{\the\headline}\vskip .5truecm \fi}
\def\makefootline{\iffrontpage\vskip  0.truecm\line{\the\footline}
               \vskip -.15truecm\line{\the\date\hfil}
              \else\line{\the\footline}\fi}
\paperheadline={
\iffrontpage\hfil
               \else
               \tenrm\hss $-$\ \folio\ $-$\hss\fi    }
\def\sectionfont{\bf}    }
%
%
%
%
\def\SIZE{\hsize=15.73truecm\vsize=23.11truecm}
\def\OFFSET{\voffset=0.0truecm\hoffset=0.truecm}
\message{DEFAULT FORMAT}
\def\papersize{\SIZE\OFFSET\skip\footins=\bigskipamount
\normaldisplayskip= 35pt plus 3pt minus 7pt}
\Pubnum={\rm NIKHEF/\the\pubnum }
\def\title#1{\vskip\frontpageskip\vskip .50truein
     \titlestyle{\seventeencp #1} \vskip\headskip\vskip\frontpageskip
     \vskip .2truein}
\def\author#1{\vskip .27truein\titlestyle{#1}\nobreak}

\def\p@bblock{\begingroup \tabskip=\hsize minus \hsize
   \baselineskip=1.5\ht\strutbox \topspace-2\baselineskip
   \halign to\hsize{\strut ##\hfil\tabskip=0pt\crcr
   \the \Pubnum\cr}\endgroup}
\def\makefootline{\iffrontpage\vskip .27truein\line{\the\footline}
                 \vskip -.1truein
              \else\line{\the\footline}\fi}
\paperfootline={\iffrontpage\message{FOOTLINE}
\hfil\else\hfil\fi}
\paperheadline={
\iffrontpage\hfil
               \else
               \twelverm\hss $-$\ \folio\ $-$\hss\fi}
%
%
\def\nup#1({\refbreak\ Nucl.\ Phys.\ $\underline {B#1}$\ (}
\def\plt#1({\refbreak\ Phys.\ Lett.\ $\underline  {#1}$\ (}
\def\cmp#1({\refbreak\ Commun.\ Math.\ Phys.\ $\underline  {#1}$\ (}
\def\prp#1({\refbreak\ Physics\ Reports\ $\underline  {#1}$\ (}
\def\prl#1({\refbreak\ Phys.\ Rev.\ Lett.\ $\underline  {#1}$\ (}
\def\prv#1({\refbreak\ Phys.\ Rev. $\underline  {D#1}$\ (}
\def\und#1({            $\underline  {#1}$\ (}
%
%

\def\rB{\hfil\penalty1000\hfilneg}
%
%
\hyphenation{sym-met-ric anti-sym-me-tric re-pa-ra-me-tri-za-tion
Lo-rentz-ian a-no-ma-ly di-men-sio-nal two-di-men-sio-nal}
%
%
%
%

\def\coeff#1#2{{\textstyle { #1 \over #2}}\displaystyle}
\def\boxit#1{\vbox{\hrule\hbox{\vrule\kern3pt
\vbox{\kern3pt#1\kern3pt}\kern3pt\vrule}\hrule}}
\message{ by V.K, W.L and A.S}
\catcode`@=12
\paperstyle

\input tables
\def\galex{\rrr\galex{J. Fuchs, A.N. Schellekens and C. Schweigert,
 Nucl. Phys. B437 (1995) 667.}}
\def\DVVV{\rrr\DVVV{R.~Dijkgraaf, C.~Vafa, E.~Verlinde and
H.~Verlinde,
\cmp 123 (1989) 16.}}
\def\CIZ {\rrr\CIZ {
A.~Cappelli, C.~Itzykson     and J.-B.~Zuber,
\nup280 (1987)  445;\rB \cmp113 (1987) 1.}}
\def\DGM {\rrr\DGM  {L.~Dolan, P.~Goddard and P.~Montague,
\plt B236 (1990) 165.}}
\def\SchM{\rrr\SchM{A.N.~Schellekens, \cmp 153 (1993)  159.}}
\def\FSSs{\rrr\FSSs{
J. Fuchs, A.N. Schellekens and C. Schweigert,
\nup473 (1996) 323.}}
\def\Gan{\rrr\Gan{T.~Gannon, Commun. Math. Phys. $\underline{161}$(1994){233};
\rB Annales Poincare Phys.Theor.65 (1996) 15-56.}}
\def\GanB{\rrr\GanB{T. Gannon, {\it The Classification 
of $SU(N)_k$ automorphism invariants},\rB hep-th/9408119;\hfill\break
T. Gannon, P. Ruelle and M.A. Walton, {\it Automorphism Modular Invariants of Current Algebras}, hep-th/9503141.}}
\def\GanK{\rrr\GanK{T. Gannon and Q. Ho-Kim, Int. J. Mod. Phys. A9 (1994) 2667.}}
\def\Mod{\rrr\Mod{T. Gannon and  M. A. Walton, {\it
Heterotic Modular Invariants and Level--Rank Duality}, 
hep-th/9804040;\hfill\break
 M. R. Abolhassani and F. Ardalan, Int. J. Mod. Phys. A9 2707 (1994);
\hfill\break
 J. Fuchs, B. Gato-Rivera, A.N. Schellekens and C. Schweigert,
\plt B334 (1994) 113;\hfill\break  A. Font  Mod.Phys.Lett.A6 (1991) 3265.
\hfill\break For older references see: 
J. Fuchs, {\it Affine Lie Algebras and Quantum Groups}
\rB
           [Monographs on Mathematical Physics],
       Cambridge University Press (1992).}}
\def\ScYe{\rrr\ScYe{A.N.~Schellekens and S.~Yankielowicz,
\nup334 (1990) 67.}}
\def\ScYE{\rrr\ScYE{A.N.~Schellekens and S.~Yankielowicz,
\nup334 (1990) 67; \rB Int.~J.~Mod.~ Phys.~\und{A5} (1990) 2903.}}
\def\DixGin{\rrr\DixGin{L. Dixon, cited in P.~Ginsparg, {\it Applied Conformal Field
Theory}, 
in the {\it Proceedings of the Les Houches school on
 Fields, Strings and Critical Phenomena, 1988.}}}
\def\CoSu{\rrr\CoSu{F.~Bais and
P.~Bouwknegt,
\nup279 (1987) 561;\hfill\break  A.N.~Schellekens and N.P.~Warner,
\prv34 (1986) 3092.}}

\def\Zbf{{\bf Z}}
\def\Qbf{{\bf Q}}
\def\mod{{\rm~mod~}}
\def\char{{\cal X}}
\pubnum={{}}
\Pubnum={{}}
\rightline{NIKHEF/98-020}
\rightline{hep-th/9806162}
\rightline{June 1998}
\date{June 1998}
\pubtype{CRAP}
\titlepage
\message{TITLE}

\title{\fourteenbf Cloning SO(N) level 2}
\author{A. N. Schellekens\foot{t58@attila.nikhef.nl}}
\line{\hfil \it NIKHEF-H, P.O. Box 41882, 1009 DB Amsterdam,
The Netherlands  \hfil}
\bigskip

\abstract \noindent

For each $N$ an infinite number of Conformal Field Theories is
presented that has the same fusion rules as $SO(N)$ level 2.
These new theories are obtained as 
extensions of the chiral algebra of $SO(NM^2)$ level 2, and 
correspond to new modular invariant partition functions of these
theories. A one-to-one map between the
$c=1$ orbifolds of radius $R^2=2r$
and $D_r$ level 2 plays an essential role.

\baselineskip= 15.0pt plus .2pt minus .1pt

\chapter{Introduction}

Since the discovery of the ADE classification of modular invariant
partition functions (MIPF's) of affine $SU(2)$ \CIZ, the generalization to other
affine Lie algebras has continued to fascinate a select group of people 
(see \eg\ \Mod\galex\SchM). 
Although completeness proofs were given in a few other cases 
\multref\Gan{\GanB\GanK}, the 
natural classification problem for affine algebras based on simple
Lie algebras remains unsolved, and the generalization to semi-simple
algebras looks totally impossible. In the ``simple" case there
may at  some point have
a belief that most MIPF's had been found,
but that belief was shaken
several times in the last few years. In particular in \galex\ an
infinite series of new MIPF's of automorphism type was described, for
the affine algebras $D_r$ and $B_r$ at level 2.

In this paper I will describe another infinite series occurring
for the same algebras, this time of extension type, and apparently not
yet known (these new extension invariants should not be confused with 
another series mentioned in \galex, which was shown to be unphysical;
the new theories we describe here are definitely physical).
The starting point of the analysis is a new interpretation of the
partition functions of \galex. 

It turns out that both types of MIPF (those of \galex\ and the
new ones described here) have a natural interpretation in terms
of $c=1$ orbifolds at arbitrary rational radius. The link between the
orbifolds and the affine algebra goes via the coset description
of the former. The $c=1$ orbifolds
of a circle are described by the coset CFT's \DixGin
$$ {SO(N)_{1} \times SO(N)_{1}\over SO(N)_{2}} \ , \eqn\coset $$
where the subscript denotes the level.
It turns out that for even $N=2r$ one obtains all radii $R^2=2r$,
whereas for odd $N=2r+1$ one only gets the radii $R^2=2(2r+1)$. In
the latter case the field identification has fixed points \ScYe, so
we will focus first on the $D_r$ affine algebras, 
where the analysis is easier. 

We introduce the 
following notation for $D_r$ Lie-algebra representations in terms of
Dynkin labels (the last two of which are the spinor labels).
$$ \eqalign{
&(0):\quad  (0,\ldots\ldots\ldots\ldots\ldots,0,0) \cr 
&(v):\quad  (1,\ldots\ldots\ldots\ldots\ldots,0,0) \cr 
&(s):\quad  (0,\ldots\ldots\ldots\ldots\ldots,0,1) \cr 
&(c):\quad  (0,\ldots\ldots\ldots\ldots\ldots,1,0) \cr 
&(\ell):\quad (0,\ldots,~0,1,0,\ldots,~0,0) \ .\cr} $$ 
where in the last line the $1$ is in position $\ell$. Furthermore
we denote by $(vv), (vs), \ldots$ \etc. the sum of the corresponding
Dynkin labels, \eg\ $(vv)=(2,0,\ldots,0,0)$, \etc. The algebra $D_{r,1}$ has
four unitary highest weight representations with ground state labels
$(0),(v),(s),(c)$, while for $D_{r,2}$ one gets in addition the
representations
 $(vv),(vc),(vs),$ $(ss),(cc),(sc)$ and $(\ell), \ell=2,r-2$. Hence 
the total number of primary fields of the $D_{r,2}$ theory is $r+7$. 
The coset theory has a four-fold field identification with
identification currents \ScYe\ $(s,s;ss)$, $(c,c;cc)$ and $(v,v;vv)$. These
currents act without fixed points, and the total number of primaries in
the coset theory is therefore $(4\times 4 \times( r+7) / 4 / 4)=r+7$.
This is the same number of fields as for the orbifolds of the $R^2=2r$ 
circle theory \DVVV, and looking at the spectra one concludes  
that the two theories are indeed the same.

The identification with the orbifold spectrum is as follows
for odd rank (column 1) and even rank (column 2).

\vskip 1.truecm
\begintable
Coset reps., $r$ odd | Coset reps., $r$ even | Orbifold reps. | Conformal weight \cr 
(0,0;0) |(0,0;0) |$[0]$ | 0 \nr
(0,0;vv) | (0,0;vv)|$[V]$ | 1 \nr
(0,v;ss) |(0,0;ss) | $[S]$ | $r/4$ \nr
(0,v;cc) |(0,0;cc) | $[C]$ | $r/4$ \nr
(0,c;s) |(0,s;s) | $ [\sigma]  $ | ${1\over16}$  \nr
(0,s;c) |(0,c;c) | $ [\tilde\sigma]  $| ${1\over16}$\nr
(0,c;vc) |(0,s;vc) | $ [\sigma']  $| ${9\over16}$ \nr
(0,s;vs) |(0,c;vs) |$ [\tilde\sigma']$| ${9\over16}$ \nr
$(0,0;\ell)$ | $(0,0;\ell)$|  $[\ell], \ell$ even | ${\ell^2\over 4r}$ \nr
$(0,v;\ell)$ |$(0,v;\ell)$ | $[\ell], \ell$ odd| ${\ell^2\over 4r}$ \endtable
\vskip 1.truecm
The representations
in columns 1 and 2 are identification orbit representatives chosen
by requiring the first entry to be $(0)$. The second entry is then fixed
by the selection rules of the coset embedding, given the third entry.
The last two columns fix our notation for the orbifold fields, whose
weight is listed in the last column.  

This table not only implies an equivalence between the coset CFT 
\coset\ and the $c=1$ orbifolds,
 but it also implies an isomorphism between the fusion rings
of the $c=1$, $R^2=2r$ orbifold and the affine algebra $D_{r,2}$. 
This isomorphism is a consequence of the fact that $D_{r,1}$ has simple
fusion and that $D_r$ fusion (at any level) preserves conjugacy class charges.
The modular transformations and fusion rules of the $c=1$ orbifolds have
been obtained in \DVVV, but I do not know if the isomorphism with
$D_{r,2}$ was noticed before.

The modular transformation matrices $S$ and $T$ of $D_{r,2}$ and the
orbifold are not the same, but differ by a complex conjugation (since
$D_{r,2}$ appears in the denominator in \coset) and some phases originating
from $D_{r,1}$. Nevertheless the connection is close enough to gain
information about modular invariant partition functions in one case and
use it in the other case.  

Non-trivial MIPF's of the $c=1$ orbifold theory can be expected to 
exist because the diagonal invariant only describes $R^2=2r$ orbifolds,
whereas the orbifold theory exists for any $R$. It was pointed
out in \DVVV\ that the chiral algebra for the orbifold $R^2=2r/q$
(with $r$ and $q$ relative prime) is the same as for $R^2=2rq$.
This implies that
the former theory must be described by a non-diagonal 
modular invariant combination of the $R^2=2rq$ characters.  
This invariant is however not discussed in \DVVV.

Let us first review the partly analogous situation 
for circle compactifications of a single real boson.
Geometrically, such compactifications are described in terms of one
modulus, the radius of the circle $R$. The description in terms of
conformal field theory is more complicated, and one has to 
distinguish three cases: (a) $R^2 \in 2\Zbf$, (b) $R^2\in 2\Qbf$ and
(c) $R^2$ irrational. In case (a) the description is in terms
of   
a chiral algebra generated by the operators ``$\partial X$" and
vertex operators related to the even lattice $\ell R^2$, $\ell \in \Zbf$.
The circle compactification corresponds to the diagonal invariant
built out of the primary fields of this extended algebra. In case (c)
the chiral algebra is just the free boson algebra $\partial X$ and there
is an infinite number of primaries. The other rational circles
(case (b)) do not correspond to diagonal invariants of some chiral
algebra (case (a) already exhausts all possible chiral algebras), 
but to automorphism invariants generated by simple currents. 

To fix the notation, define ${\cal U}_{2r}$ as the $U(1)$ CFT with
$2r$ primaries, which has as its extended algebra the one
corresponding to $R^2=2r$ above. All the primaries are simple
currents, which we will label by $J$, $J=0,\ldots,2r-1$.
Suppose $r=pq$, where $p$ and
$q$ are relative prime.
Now consider the modular invariant partition function
generated by the simple current $J=2p$; this is an automorphism
invariant.
It is easy to show that
this yields precisely the circle compactification with radius 
$R^2=2p/q$. The current $J=2q$ gives radius $R^2=2q/p$, the
T-dual ($R \leftrightarrow 2/R$) of the previous case. In particular
$J=2$ yields the T-dual of the $R^2=2r$ circle compactification, 
and corresponds on the other hand to the charge conjugation invariant
of ${\cal U}_{2r}$. 

This exhausts the set of rational circles, but not the set of simple
currents. One may in fact use any simple current $J=2j$ (odd $J$'s are
not in the effective center \ScYE), and in general one obtains $R^2=2r/N^2$,
where $N$ is the order of $J$. If $r$ contains $N^2$ as a factor this
operation reduces the theory to one with $R^2=2\tilde r$,
$\tilde r=r/N^2 \in \Zbf$. This should therefore be a diagonal theory.
Indeed, in that case the current $J$ has integer spin and generates an
extension of the chiral algebra (thus reducing $r$ to $\tilde r$),
and not an automorphism (there are also mixed cases where a factor
$\tilde N^2$ of $N^2$ divides $r$).

One would expect a similar situation to occur in the orbifold 
theory, but there is a clear difference: the orbifold has
only four simple currents, the fields $[0],[V],[S]$ and $[C]$ (except
for $r=1$, when the orbifold coincides with the $R^2=8$ circle). This
is clearly not sufficient, and we are forced to conclude that at least
the invariants required to go from $R^2=2pq$ to $R^2=2p/q$ must be 
exceptional (\ie\ not simple current) invariants. 
Because of the isomorphism, their existence
then implies the existence
of analogous invariants for $D_{r,2}$.     

At this point we make contact with \galex, where precisely these
invariants were discovered from an entirely different point of view,
namely Galois symmetry. They exist for $D_r$ and also for $B_r$. 
The invariants described in \galex\ occur whenever 
$2r-1=\Pi_{i=1}^K p_i$ (for $B_r$)
or $r=\Pi_{i=1}^K p_i$ (for $D_r$), where the factors $p_i$ are
distinct primes. In total there are $2^{K-1}$ distinct
invariants, including
the identity, precisely one for each way of writing $2r-1$ (resp. $r$)
as a product of two factors. It is easy to check that these
automorphisms do indeed yield the partition function one expects for
the $R^2=2p/q$ orbifolds. On closer inspection the results of \galex\
can be generalized to all cases where $r=pq$ (or $2r-1=pq$) if $p$ and
$q$ are relative prime.  This exhausts\foot{The $c=1$ theories, just as
$D_r$ level 2, have (at least) one additional automorphism, corresponding
in $D_r$ to spinor conjugation (or charge conjugation for $r$ odd). This
automorphism invariant may be interpreted as the $T$-dual of the orbifold
theory. For $r=1$, the three-critical point,  
it becomes the usual $T$-duality of the circle.}all rational orbifolds, and 
gives us a {\it raison d'etre} of these exceptional $D_r$
automorphisms.

In addition to these automorphisms the circle theories at suitable radii
also have extensions by integer spin simple currents. Although nothing
requires their orbifold equivalent to exist (we already found all
rational orbifolds) it turns out that they {\it do} nevertheless exist.
They are given by the following  
general formula for the $D_{r,2}$ invariant, if
$r=\tilde r M^2$ and $M$ is odd.
$$\eqalign{ &\mid\char_{0}  +\sum_{m=1}^{(M-1)/2} \char_{2m\tilde r M}\mid^2
+\mid \char_{vv}  +\sum_{m=1}^{(M-1)/2} \char_{2m\tilde r M}\mid^2 \cr
           +&\mid\char_{ss}+\sum_{m=1}^{(M-1)/2} \char_{(2m-1)\tilde r M}\mid^2+
         \mid\char_{cc}+\sum_{m=1}^{(M-1)/2} \char_{(2m-1)\tilde r M}\mid^2\cr
&+ \sum_{l=1}^{\tilde r-1} 
\mid \sum_{m=0}^{M-1} \char_{r-|r-lM-2m\tilde r M|}\mid^2\cr
+&\mid\char_{s}\mid^2
 +\mid\char_{c}\mid^2
 +\mid\char_{vc}\mid^2
 +\mid\char_{vs}\mid^2\cr}\eqn\Dinv $$
There is a completely analogous invariant for the $c=1$ orbifold, but of
course it merely describes the $c=1$ orbifold at a reduced radius, which
is nothing new. 
The $D_{r,2}$ invariant, however, describes a conformal field theory 
with fusion rules that are isomorphic to $D_{\tilde r,2}$, but with 
a {\it different} spectrum. Hence for any $\tilde r$ we get 
an infinite series of 
CFT's that are ``fusion-isomorphic"
to $D_{\tilde r,2}$. These new CFT's are labelled by an odd integer $M$, 
and have a chiral algebra that
is an extension of $D_{\tilde rM^2,2}$. We will call this
new theory $D^M_{\tilde r,2}$.  

To check that the proposed partition function is indeed modular
invariant we need the transformation matrix $S$. It can be computed 
either using the Kac-Peterson formula for the affine algebras, but 
can be more easily obtained from the results of \DVVV, adding a few
$D_{r,1}$ phases where needed. 
Checking the invariance is then straightforward. Note that
the partition function has a somewhat unusual form: the characters of the
extended theory are linear combinations of either $\half(M+1)$, $M$ or
one character of the original theory, and furthermore the first four
blocks are not ``orthogonal", \eg\ the currents that extend the algebra
occur also in a non-identity character. Furthermore the short blocks
have multiplicity 1, unlike fixed point characters of simple current
invariants. The fact that all multiplicities are 1 makes it straightforward
to compute the matrix $S$ of the extended theory. As expected, for the 
orbifold theory it is again the one given in \DVVV, whereas for 
$D_{r,2}$ it differs by the $D_{r,1}$ phases. Since the matrices 
$S$ for $D_{\tilde rM^2,1}$ and $D_{\tilde r}$ are identical if $M$
is odd (because $M^2=1 \mod 8$), we may replace the $D_{r,1}$  by 
$D_{\tilde r,1}$ phases, which shows that $D^M_{\tilde r,2}$ has
the same matrix $S$ as $D^1_{\tilde r,2}=D_{\tilde r,2}$.   

The foregoing results can be extended from odd $M$ to all integer $M$.
Consider first $M=2$. In this case
one does not require an exceptional
invariant, but instead the simple current $[S]$ (or $[C]$) achieves
the reduction of $r$ by a factor 4. The invariant generated by
this simple current is\foot{The invariant is written here
 for $D_{r,2}$, but the corresponding
invariant for the $R^2=2r$ orbifold is completely analogous.}
$$ \eqalign{ &\mid \char_{0} + \char_{ss}\mid^2 + 
\mid \char_{vv} + \char_{cc}\mid^2\cr
&+ 2 \mid \char_s^2 \mid^2 + 2 \mid \char_{vc} \mid^2 
+ 2 \mid \char_{r/2} \mid^2 \cr &+ 
\sum_{l=2\atop \ell even}^{r/2-2}
 \mid \char_{\ell} + \char_{r-\ell}\mid^2\cr} \eqn\scinv $$
Here each fixed point splits into two distinct fields, and
hence we get a total of $2+6+(r/4-1)=\tilde r+7$ fields in the new theory.
In this case a fixed point resolution procedure is needed to compute
the new matrix $S$, but the required formalism is available \ScYe\FSSs. 
Fixed point resolution then   
yields the matrix $S$ of $D_{\tilde r,2}$, $\tilde r = r/4$ (up to
a few phase changes if $\tilde r \not= 0 \mod 4$ 
due to the fact that the matrices $S$ for $D_{r,1}$
and $D_{\tilde r,1}$ are not identical; there are no such phase changes
in the orbifold case, because then \scinv\ produces a theory that is
{\it identical} to the $R^2=2\tilde r$ orbifold). 

Since we always land on another $D_{r,2}$ theory these extensions
can be performed successively, which allows a reduction of $r$
by any factor $4^n$. A reduction by a factor 16 or more is however
not a standard simple current invariant. It turns out that the first
reduction by a factor 4 promotes the primary field $(r/2)$
 to a simple current,
which is then used in the second stage. This is possible because $(r/2)$
is a fixed point in the first stage. Finally we note that when removing
the last factor of 4 from $r$ we pass from the $r$ even to $r$ odd.
In \DVVV\ it was found that $S$ is real for $r$ even, but complex for
$r$ odd. The imaginary part comes from the fixed point resolution.
Without fixed points the new matrix
elements of $S$ are real linear combinations of the old ones, and
an imaginary part cannot be generated. Using the table in \FSSs\ one
readily finds that an imaginary fixed point resolution matrix occurs
precisely for $D_{r}$, $r=0\mod 4$, $r\not=0\mod 8$.      

\section{The series for $B_{r,2}$}

For the algebra $B_{r,2}$ we can also derive the matrix $S$ from the
results of \DVVV.   
In this case we use the coset theory
$$ {B_{s,1} \times B_{s,1}\over B_{s,2}} \eqn\coset $$
with identification current $(v,v;vv)$. The identification of the 
spectrum is now slightly more difficult because this current has a fixed
point, but is not difficult to show that this coset gives the $c=1$
orbifold with $r=2(2s+1)$. Hence using $B_{s}$
 (\ie\ $SO(L), L = 2s+1$) we only get
 cosets with even radii $r=2L$. 
The $B_{s,2}$ spectrum consists of the representations
$(0),(vv),(s),(sv),(\ell)$, $\ell=1,\ldots,(L-1)/2$, where we use
a notation analogous to $D_r$. The representation $(ss)$ is in fact an
The anti-symmetric tensor of rank $(L-1)/2$ and is more conveniently
denoted as such.     
The identification
of the spectrum is \vskip 1.truecm
\begintable
Coset repr. | $c=1$ repr. |Coset repr. | $c=1$ repr.\cr
(0,0;0) | [0]  | $(0,0;\ell)$ | $[2\ell]$, $\ell$ even \nr
(0,0;vv) | $[V]$  | $(0,0;\ell)$ | $[2L-2\ell]$, $\ell$ odd \nr
(0,v;0) | $    [S] $ |  $(0,v;\ell)$ | $[2L-2\ell]$, $\ell$ even\nr
(v,0;0) | $    [C] $ |  $(0,v;\ell)$ | $[2\ell]$, $\ell$ odd\nr
(0,s;s) | $ [\sigma]  $  |  $(s,s;0)$  | $[L]$ \nr

(s,0;s) | $ [\tilde\sigma]  $ | $(s,s;\ell)_1$ | $[L-2\ell]$ \nr
(0,s;vs) | $ [\sigma']  $ |  $(s,s;\ell)_2$  | $[L+2\ell]$ \nr
(s,0;vs) |$ [\tilde\sigma']$ | | \endtable
\vskip 1.truecm
Choosing a basis $(0,vv,\ell,s,vs)$ for the $B_{s,2}$ representations 
we get then the following result for $S$
$$ S=\pmatrix{ a & a & 2a & \half & \half \cr
             a & a & 2a & -\half & -\half \cr
             2a & 2a & 4a \cos{2\pi\ell\ell'\over L} & 0 & 0 \cr
             \half & -\half & 0  & \half & -\half \cr
             \half & -\half & 0  & -\half & \half \cr} \eqn \Bmat$$
where $a={1\over 2 \sqrt L}$. 
Now consider $L=\tilde L M^2$. In that case the following is a modular
invariant partition function
$$\eqalign{ &\mid \char_0 + \sum_{m=1}^{(M-1)/2} \char_{m L M} \mid^2
+\mid \char_{vv} + \sum_{m=1}^{(M-1)/2} \char_{m L M} \mid^2\cr
&+ \mid \char_{s}\mid^2+\mid \char_{sv}\mid^2\cr
&+ \sum_{l=1}^{(\tilde L-1)/2} 
\mid \sum_{m=0}^{M-1} \char_{\half L-|\half L-lM-m\tilde L M|}\mid^2\cr}$$
The subscript of the last term may be simplified by doubling the range
of the anti-symmetric tensors from $1 \leq l \leq (L-1)/2$ to 
$1 \leq l \leq (L-1)$, 
identifying $l$ with $L-l$ (note that the formula for $S$ is
invariant under this). Then the subscript is simply $lM+m\tilde L M$
(a similar remark applies to \Dinv). 

The matrix $S$ that transforms the new, extended characters turns out to be
precisely the one of $SO(\tilde L)$. In the
limiting case $\tilde L=1$ this requires an extension to ``$SO(1)$", but
\Bmat\ is well-defined in that limit. In the limit one gets a 
$D_{r,1}$-type matrix, and this implies that all four representations
are simple currents.  
In particular $(s)$ and $(sv)$ become simple
currents after the extension of $SO(M^2)$. The conformal weights of 
these two representations are $\coeff1{16}(M^2-1)$ and 
$\coeff1{16}(M^2+7)$ respectively, and for any odd $M$ one of these
weights is integer and the other half-integer. The integer spin simple
current can be used to extended the algebra even further, leaving
only a single representation. The result is a meromorphic CFT \DGM\SchM. Indeed,
the first examples of this phenomenon are already known, and correspond
to the conformal embedding $B_{4,2} \subset E_{8,1}$ and the $c=24$
meromorphic CFT based on $B_{12,2}$ (Also the $D_{9,2}$ extension
can be found on the list of $c=24$ meromorphic CFT's, although 
less directly).

It should be stressed that these 
new MIPF's are neither simple current
invariants nor conformal embeddings, \ie\ their chiral algebras 
contain neither simple currents nor spin-1 currents (apart from those
of the original affine Lie algebra). Nevertheless they are closely
related to both types. The new theories all contain a representation
$[V]$ which is a spin-1 simple current. They can therefore be further
extended, and then their $SO(N)_2$ algebra is extended to $SU(N)_1$. This 
additional extension is both a simple current invariant and a 
conformal embedding. The resulting theory is a non-trivial MIPF of
$SU(N)_1$, and is in fact a simple current invariant of $SU(N)_1$. 
Hence the new invariants are related to already known ones in the 
following way 
$$\eqalign{ SU(N)_1\quad &{{\lower10.pt\hbox{$\longrightarrow$}} \atop {\scriptstyle S.C }}\quad SU(N)_1^{\rm ext} \cr
    \bigg\uparrow \scriptstyle C.E\quad\quad &\quad\quad\quad\quad
 \bigg\downarrow \cr
SO(N)_2\quad &{{\lower10.pt\hbox
{$\longrightarrow$}} \atop \scriptstyle H.S.E }\quad SO(N)_2^{\rm ext}\cr}
\eqn\vdiag $$
Here ``S.C.", ``C.E" and ``H.S.E" stand for simple current,
conformal embedding and higher spin extension respectively. The new
MIPF's are in the lower right corner.
The inverse of the conformal embedding is a $\Zbf_2$ orbifold (inversion
of the Cartan sub-algebra of $SU(N)_1$).
 Applying this same orbifold
procedure to the simple current extensions corresponds to the fourth,
unmarked arrow in the diagram: the new $SO(N)_2$ invariants are
$\Zbf_2$ orbifolds of simple current extensions of $SU(N)_1$. This 
proves that a sensible CFT corresponding to these MIPF's exists.  
The $\Zbf_2$-orbifold was in fact used in \DGM\ to construct 
-- among others -- the meromorphic $B_{12,2}$ theory, described above,
from the $A_{24}$ self-dual lattice. 

The diagram \vdiag\ suggests a generalization to other cases
of conformal subalgebras $H\subset G$, where $G$ can be extended by
simple currents. 
The non-trivial issue is the existence of the unmarked arrow in the
diagram. This corresponds to ``undoing" the conformal embedding for
the simple current extension of $G$. In the example discussed in this paper,
this undoing is achieved by a simple and well-understood orbifold
procedure, but this is not true in general. Furthermore $H$ has
representations not present in $G$ (corresponding to twisted states in
an orbifold), and the matrix $S$ of $G$ gives no information
about $S$ on these states. A partial inspection of the list of
conformal embeddings \CoSu\ indicates that \vdiag\ does not generalize,
or at least not to all cases.

\ack

I would like to thank
Christoph Schweigert for comments and for reading the manu\-script.

\par \penalty-4000\vskip\chapterskip
   \spacecheck\referenceminspace \immediate\closeout\referencewrite
   \referenceopenfalse
   \line{\fourteenrm\hfil REFERENCES\hfil}\vskip\headskip
   \endlinechar=-1
   \input referenc.texauxil
   \endlinechar=13
   
\end


Consider now the $\Zbf_2$ orbifold of the circle. 
Again the geometric description has the advantage of being
continuous in $R$, but the CFT description is instructive as well.
For $R^2=2P$
the orbifolds are described as diagonal invariants of conformal
field theories with $P+7$ primaries (the vacuum, the
field "$\partial X$" with weight one,two twist fields with
weight $1/16$, two exited twist fields with weight $9/16$, two
fields with weight $P/4$ from the untwisted sector of the $J=P$ circle
representation, and $P-1$ fields from the untwisted sector
of the circle fields with $0 < |J| < P$.). On general grounds one
expects that orbifolds with arbitrary rational radius also have
a CFT description, but in contrast to the circle theories, the
orbifold theories do not have simple currents that can achieve this
(the only simple current is $\partial X$, but using this just
extends the chiral algebra to the circle theory).

The solution to this puzzle is that the orbifold CFT's for 
certain values of $R^2$ have exceptional automorphism invariants.
This fact turns out to be closely related to the $SO(N)$ invariants
discussed earlier, and follows from the coset description of the
orbifold theories. The coset theory is
$$ {SO(N)_1 \times SO(N)_1 \over SO(N)_2} $$
For $N$ odd the determination of the spectrum of this coset CFT
requires a two-fold field identification with a total of
$\half(N-1)$ fixed points.  The total number of primary fields turns
out to be $2N+7$, and the spectrum matches the orbifold theory 
for $P=2N$. For $N$ even there is four-fold identification without
fixed points, and there are $N/2+7$ primaries, matching the orbifold
with $P=N/2$. Note that for even $P$ the orbifold is described both
by $N$ odd and $N$ even. 

Having written the orbifold theory in terms of $SO(N)_2$ we can
immediately make contact with the automorphism invariants of 
those theories discussed earlier. Consider first $N$ odd, 
with $N$ a multiple of $K$ single primes. 
Then $R^2=2P=4N$.
There are 
$2^{K-1}$ invariants, corresponding to the number of ways
of writing $N=pq$. These invariants correspond precisely
to the orbifold theories with $R^2=2p/q$ (or equivalently
$R^2=2q/p$; for the $c=1$ orbifolds
T-duality is trivial; the theories with  $R^2=2q/p$
$R^2=2p/q$ are not just dual but in fact identical. This is due
to the fact that T-duality for the circle theory is charge 
conjugation, and charge conjugation projects to the identity in the
orbifold theory.)

Now consider $N$ even, also a product of $K$ single primes.

If $r=4\tilde r $ the reduction can be achieved by a simple current.
The two spinor currents $(ss)$ and $(cc)$ then have integer spin $r/4$,
and each of them can be used to extend the chiral algebra. This
simple current invariant has the following form
$$ \eqalign{ &\mid \char_{0} + \char_{ss}\mid^2 + 
\mid \char_{vv} + \char_{cc}\mid^2\cr
&+ 2 \mid \char_s^2 \mid^2 + 2 \mid \char_{vc} \mid^2 
+ 2 \mid \char_{r/2} \mid^2 \cr &+ 
\sum_{\ell \hbox{~even}; l=2}^{r/2-2}
 \mid \char_{\ell} + \char_{r-\ell}\mid^2\cr} $$
In this case each fixed point splits into two distinct fields, and
hence we get a total of $2+6+(r/4-1)=\tilde r+7$ fields in the new theory.

This then corresponds to the reduction to a theory that is isomorphic to
$D_{\tilde r}$. After the extension the fixed point field 
$\char_{r/2}$ becomes a simple current. Two in fact, corresponding to the
fields $[S]$ and $[C]$ of the new theory. 

These invariants (to be precise, all the $B_n$ and half of the $D_n$
invariants) were discovered in \galex\ using Galois symmetry, but it was
pointed out that they could also be understood in terms of the 
conformal embedding $SO(N)_2 \subset SU(N)_1$, where the subscript
denotes the level.